\documentclass[letterpaper, 10 pt, conference]{ieeeconf}
\IEEEoverridecommandlockouts
\usepackage{amsmath}
\usepackage{amssymb}
\usepackage{epsfig}

\title{\LARGE \bf Motion camouflage in three dimensions}

\author{P.~V. Reddy, E.~W. Justh, and P.~S. Krishnaprasad
\thanks{
This research was supported in part by the Naval Research Laboratory under
Grants No.~N00173-02-1G002, N00173-03-1G001, N00173-03-1G019, and
N00173-04-1G014; by the
Air Force Office of Scientific Research under AFOSR Grants
No.~F49620-01-0415 and FA95500410130; by the Army Research Office
under ODDR\&E MURI01 Program Grant No.~DAAD19-01-1-0465 to the Center for
Communicating Networked Control Systems (through Boston University);
and by NIH-NIBIB grant
1 R01 EB004750-01, as part of the NSF/NIH Collaborative
Research in Computational Neuroscience Program.}
\thanks{P.V. Reddy and P.S. Krishnaprasad are with the Institute for
Systems Research and the Department
of Electrical and Computer Engineering at the University of
Maryland, College Park, MD 20742, USA. {\tt\small vishwa@umd.edu, 
krishna@umd.edu}}
\thanks{E.W. Justh is with the Institute for Systems Research at
the University of Maryland, College Park, MD 20742, USA.  
{\tt\small justh@umd.edu}}
}

\begin{document}

\maketitle
\thispagestyle{empty}
\pagestyle{empty}

\begin{abstract}

We formulate and analyze a three-dimensional model of 
motion camouflage, a stealth strategy observed in nature.
A high-gain feedback law for motion camouflage is formulated
in which the pursuer and evader trajectories are described
using natural Frenet frames (or relatively parallel adapted
frames), and the corresponding natural curvatures serve as controls.
The biological plausibility of the feedback law is discussed,
as is its connection to missile guidance.  Simulations
illustrating motion camouflage are also presented.
This paper builds on recent work on motion camouflage
in the planar setting \cite{mcarxiv05}.

\end{abstract}

\section{Introduction}

Motion camouflage is a stealth strategy employed by various visual
insects and animals to achieve prey capture, mating or territorial  
combat. In one type of motion camouflage, the predator 
camouflages itself against a fixed background object
so that the prey observes no relative motion between 
the predator and the fixed object.
In the other type of motion camouflage, the predator 
approaches the prey such that from the point
of view of the prey, the predator always appears to
be at the same bearing.  (In this case, we say that the
object against which the predator is camouflaged is
the point at infinity.)  For background on motion camouflage,
see \cite{mcarxiv05} and the references therein.  Motion camouflage
behavior in insects is described in \cite{srini95} 
(based on earlier work in \cite{collett75} on hoverflies)
and in \cite{srini03} (for dragonflies).  Related themes in 
insect vision and flight control are also found in \cite{srini04}.

The essential features of motion camouflage are not limited to
visual insects.  Recent work on the neuroethology of insect-capture
behavior in echolocating bats reveals a strategy geometrically 
indistinguishable from motion camouflage, referred to as the 
``constant absolute target direction'' (CATD) strategy \cite{ghose05}.
Because the bat under study, {\it Eptesicus fuscus}, hunts at night,
there is no reason to suppose that camouflage (i.e., misleading
its prey's visual system) is the bat's goal in using the CATD
strategy.  In this paper, we are concerned with describing in 
the simplest possible, biologically plausible way {\it how}
the motion camouflage or CATD strategy can be achieved using
feedback control.  This is a small first step toward 
understanding  the much more difficult question of {\it why} an animal 
like the bat {\it Eptesicus fuscus} uses such a strategy.

What sets this work apart is the structured approach used to 
derive feedback laws for motion control in three dimensions.
We model the pursuer (i.e., predator) and evader (i.e., prey)
as point particles subject to curvature (steering) control.
Although the speeds of the particles may vary, this variation
is considered to result primarily from flight conditions the
animal experiences - not primarily as a result of explicit speed
control for purposes of achieving motion camouflage.  Indeed, the 
feedback law we derive for motion camouflage
is well-defined for constant-speed motion.  However, for comparing
the theoretical feedback law to the experimentally-derived bat
trajectory data, it is useful to retain speed variability in the
model, since speed variations on the order of 50 percent are observed
as the bat maneuvers. 

This focus on systematic formulation and analysis of biologically
plausible feedback laws for motion camouflage is a distinguishing
feature of our work.  For example, in \cite{glendinning} 
motion camouflage trajectories are studied, but without explicitly
providing feedback laws which give rise to them.  In \cite{andersonneural},
feedback using neural networks is used to achieve motion camouflage,
but our approach has the advantage of giving an explicit form
and straightforward physical interpretation for the feedback control law.

In earlier work, motion camouflage in the planar setting was studied,
and a feedback law to achieve motion camouflage was derived \cite{mcarxiv05}.
The name given to the feedback law was 
{\it motion camouflage proportional guidance} (MCPG).
Here, we extend this work by formulating the problem in three
dimensions and generalizing the feedback law to the three dimensional
setting.  The
key is to describe the particle trajectories using natural Frenet frames 
\cite{bishop} - the same approach demonstrated successfully in the context 
of formation control for constant-speed particles \cite{cdc05}.
This formulation can also be used to describe missile guidance,
specifically, pure proportional navigation guidance (PPNG)
\cite{songha,oh_ha99,shneydor98}, cleanly in three dimensions.

\section{Pursuit-evasion model}

\begin{figure}
\hspace{1cm}
\epsfxsize=7cm
\epsfbox{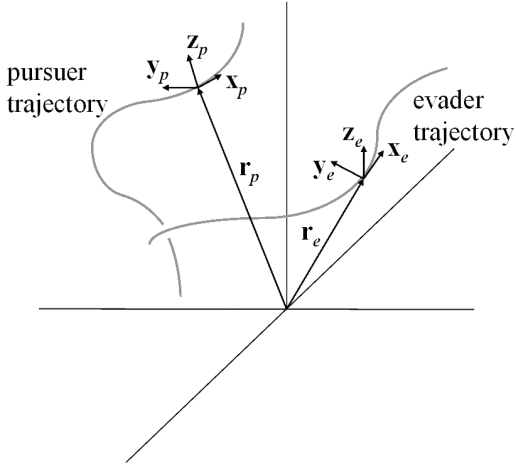}
\caption{\label{framefig3d} Trajectories for the pursuer
and evader, and their respective natural Frenet frames.
The position of the pursuer is $ {\bf r}_p $, and its
natural Frenet frame is $ \{{\bf x}_p,{\bf y}_p,{\bf z}_p \} $,
where $ {\bf x}_p $ is the unit tangent vector to its trajectory,
and $ \{ {\bf y}_p, {\bf z}_p \} $ span the corresponding normal plane
(and similarly for the evader).  The pursuer moves with speed $ \nu_p $,
and the evader with speed $ \nu_e $.}
\end{figure}

For concreteness, we consider the problem of motion
camouflage in which the predator (which we refer to 
as the ``pursuer'') attempts to intercept the prey
(which we refer to as the ``evader'') while appearing
to the prey as though it is always at the same bearing  
(i.e., motion camouflaged against the point at infinity).
The dynamics of the pursuer are given by
\begin{eqnarray}
\label{pursuer3d}
\dot{\bf r}_p \hspace{-.2cm} & = & \hspace{-.2cm}
 \nu_p {\bf x}_p, \nonumber \\
\dot{\bf x}_p \hspace{-.2cm} & = & \hspace{-.2cm}
 \nu_p ({\bf y}_p u_p + {\bf z}_p v_p),  \nonumber \\
\dot{\bf y}_p \hspace{-.2cm} & = & \hspace{-.2cm}
 -\nu_p {\bf x}_p u_p, \nonumber \\
\dot{\bf z}_p \hspace{-.2cm} & = & \hspace{-.2cm}
 -\nu_p {\bf x}_p v_p,
\end{eqnarray}
where $ {\bf r}_p $ is the position of the pursuer, $ \nu_p $ is the
speed of the pursuer,
$ {\bf x}_p $ is the unit tangent vector to the trajectory of
the pursuer, $ {\bf y}_p $ and $ {\bf z}_p $ span the normal 
plane to $ {\bf x}_p $
(completing a right-handed orthonormal basis with $ {\bf x}_p $),
and the natural curvatures $ u_p $ and $ v_p $ are the controls for the
pursuer.  Similarly, the dynamics of the evader are
\begin{eqnarray}
\label{evader3d}
\dot{\bf r}_e \hspace{-.2cm} & = & \hspace{-.2cm}
 \nu_e {\bf x}_e, \nonumber \\
\dot{\bf x}_e \hspace{-.2cm} & = & \hspace{-.2cm}
 \nu_e ({\bf y}_e u_e + {\bf z}_e v_e), \nonumber \\
\dot{\bf y}_e \hspace{-.2cm} & = & \hspace{-.2cm}
 -\nu_e {\bf x}_e u_e, \nonumber \\
\dot{\bf z}_e \hspace{-.2cm} & = & \hspace{-.2cm}
 -\nu_e {\bf x}_e v_e,
\end{eqnarray}
where $ {\bf r}_e $ is the position of the evader, $ \nu_e $ is the
speed of the evader,
$ {\bf x}_e $ is the unit tangent vector to the trajectory of 
the evader, $ {\bf y}_e $  and $ {\bf z}_e $ span the normal 
plane to $ {\bf x}_e $
(completing a right-handed orthonormal basis with $ {\bf x}_e $),
and the natural curvatures $ u_e $ and $ v_e $ are the controls for the
evader.
Figure \ref{framefig3d} illustrates equations (\ref{pursuer3d})
and (\ref{evader3d}).  Note that $ \{ {\bf x}_p, {\bf y}_p, {\bf z}_p \} $
and $ \{ {\bf x}_e, {\bf y}_e, {\bf z}_e \} $ are natural Frenet frames
(also known as relatively parallel adapted frames)
for the trajectories of the pursuer and evader, respectively \cite{bishop}.

We model the pursuer and evader as point particles, and use 
natural frames and curvature controls to describe
their motion, because this is a simple model for which we
can derive both physical intuition and concrete control laws.
Flying insects and animals (also unmanned aerial vehicles) have
limited maneuverability and must maintain sufficient airspeed to
stay aloft, so modeling them in this way is physically reasonable, 
at least for some range of flight conditions.  

Note that the forces supplied by the curvature controls 
are perpendicular to the
instantaneous direction of motion, and therefore do not
change the speed: these forces are {\it gyroscopic} forces.
However, in (\ref{pursuer3d}) and (\ref{evader3d}) we do
allow for the possibility of speed variations, as well.

\subsection{Characterizing motion camouflage}

Motion camouflage with respect to the point at infinity is given by
\cite{mcarxiv05}
\begin{equation}
{\bf r}_p = {\bf r}_e + \lambda {\bf r}_\infty,
\end{equation}
where $ {\bf r}_\infty $ is a fixed unit vector and $ \lambda $ is
a time-dependent scalar (see also Section 5 of \cite{glendinning}).

Let 
\begin{equation}
\label{rdefnplanar}
{\bf r} = {\bf r}_p - {\bf r}_e 
\end{equation}
be the vector from the
evader to the pursuer.  We refer to $ {\bf r} $ as the 
``baseline vector,'' and $ |{\bf r}| $ as the ``baseline length.''
We restrict attention to non-collision states, i.e., $ {\bf r} \ne 0 $.
In that case, the component of the pursuer velocity $ \dot{\bf r}_p $
transverse to the base line is
\begin{equation*}
\dot{\bf r}_p - \left( \frac{\bf r}{|{\bf r}|} \cdot \dot{\bf r}_p
 \right) \frac{\bf r}{|{\bf r}|},
\end{equation*}
and similarly, that of the evader is
\begin{equation*}
\dot{\bf r}_e - \left( \frac{\bf r}{|{\bf r}|} \cdot \dot{\bf r}_e
 \right) \frac{\bf r}{|{\bf r}|}.
\end{equation*}
The {\it relative} transverse component is
\begin{eqnarray}
\label{wdefn}
{\bf w} \hspace{-.2cm} & = & \hspace{-.2cm}
\left(\dot{\bf r}_p - \dot{\bf r}_e\right)
 - \left( \frac{\bf r}{|{\bf r}|} \cdot \left(\dot{\bf r}_p - \dot{\bf r}_e
 \right) \right) \frac{\bf r}{|{\bf r}|} \nonumber \\
\hspace{-.2cm} & = & \hspace{-.2cm}
\dot{\bf r} - \left( \frac{\bf r}{|{\bf r}|} \cdot \dot{\bf r}
 \right) \frac{\bf r}{|{\bf r}|}.
\end{eqnarray}

\vspace{.5cm}

\noindent
{\bf Lemma} (Infinitesimal characterization of motion camouflage):
The pursuit-evasion system (\ref{pursuer3d}), (\ref{evader3d})
is in a state of motion camouflage without collision on
an interval iff $ {\bf w} = {\bf 0} $ on that interval.

\vspace{.5cm}

\noindent
{\bf Proof}: $ (\Longrightarrow) $ Suppose motion camouflage holds.
Thus
\begin{equation}
{\bf r}(t) = \lambda(t) {\bf r}_\infty, \;\; t \in [0,T].
\end{equation}
Differentiating, $ \dot{\bf r} = \dot{\lambda} {\bf r}_\infty $.  Hence,
\begin{eqnarray}
{\bf w} \hspace{-.2cm} & = & \hspace{-.2cm}
 \dot{\bf r} - \left( \frac{\bf r}{|{\bf r}|} \cdot \dot{\bf r}
 \right) \frac{\bf r}{|{\bf r}|} \nonumber \\
\hspace{-.2cm} & = & \hspace{-.2cm}
\dot{\lambda} {\bf r}_\infty - \left( \frac{\lambda}{|\lambda|} {\bf r}_\infty
 \cdot \dot{\lambda} {\bf r}_\infty \right) \frac{\lambda}{|\lambda|}
 {\bf r}_\infty \nonumber \\
 \hspace{-.2cm} & = & \hspace{-.2cm}
 {\bf 0}  \mbox{ on } [0,T].
\end{eqnarray}
$ (\Longleftarrow) $ Suppose $ {\bf w}= {\bf 0} $ on $ [0,T] $.  Thus
\begin{equation}
\dot{\bf r} =  \left( \frac{\bf r}{|{\bf r}|} \cdot \dot{\bf r}
 \right) \frac{\bf r}{|{\bf r}|} \triangleq \xi {\bf r},
\end{equation}
so that
\begin{eqnarray}
{\bf r}(t) \hspace{-.2cm} & = & \hspace{-.2cm}
 \exp\left( \int_0^t \xi(\sigma) d\sigma \right) {\bf r}(0)
 \nonumber \\
\hspace{-.2cm} & = & \hspace{-.2cm}
 |{\bf r}(0)| \exp\left( \int_0^t \xi(\sigma) d\sigma \right)
 \frac{{\bf r}(0)}{|{\bf r}(0)|} \nonumber \\
\hspace{-.2cm} & = & \hspace{-.2cm}
 \lambda(t) {\bf r}_\infty,
\end{eqnarray}
where $ {\bf r}_\infty  = {\bf r}(0)/|{\bf r}(0)| $ and
$ \lambda(t) = |{\bf r}(0)| \exp\left( \int_0^t \xi(\sigma) d\sigma \right) $.
$ \Box $

\vspace{.5cm}

\noindent
{\bf Remark}: The above {\bf Lemma} and its proof are identical 
to the corresponding {\bf Lemma} and proof in
\cite{mcarxiv05}, but with the vectors interpreted as
three-dimensional rather than planar vectors.  $ \Box $

\vspace{.5cm}

\begin{figure}
\epsfxsize=8.5cm
\epsfbox{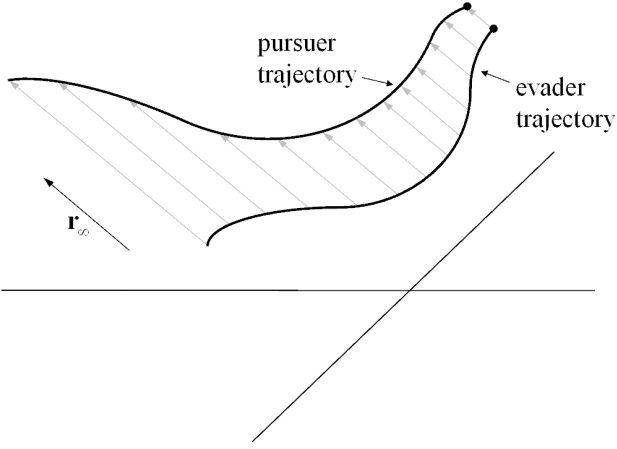}
\caption{\label{mcstatefig} Pursuer and evader trajectories in a
state of motion camouflage with respect to the point at
infinity, i.e., satisfying $ {\bf r}_p - {\bf r}_e = \lambda {\bf r}_\infty $,
where $ r_{\infty} $ is fixed and $ \lambda $ varies with time.  The
light gray vectors are baseline vectors at different instants of time:
note that they are all parallel to one another.}
\end{figure}

Figure \ref{mcstatefig} illustrates the pursuer and evader in 
a state of motion camouflage with respect to the point at infinity.

\subsection{Measuring departure from motion camouflage}

Consider the ratio
\begin{equation}
\label{gammadefn}
\Gamma(t) = \frac{ \frac{d}{dt}{|{\bf r}|}}{\left|\frac{d{\bf r}}{dt}\right|},
\end{equation}
which compares the rate of change of the
baseline length to the absolute rate of change of the baseline vector
\cite{mcarxiv05}.
If the baseline experiences pure lengthening, then the ratio assumes
its maximum value, $ \Gamma(t) = 1 $.
If the baseline experiences pure shortening, then the ratio assumes
its minimum value, $ \Gamma(t) = -1 $.  If the baseline experiences pure
rotation, but remains the same length, then $ \Gamma(t) = 0 $.
Noting that
\begin{equation}
\label{ddtnormr}
\frac{d}{dt}{|{\bf r}|} = \frac{\bf r}{|{\bf r}|}\cdot \dot{\bf r},
\end{equation}
we see that $ \Gamma(t) $ may alternatively be written as
\begin{equation}
\label{gammadotprod}
\Gamma(t) = \frac{\bf r}{|{\bf r}|} \cdot \frac{\dot{\bf r}}{|\dot{\bf r}|}.
\end{equation}
Thus, $ \Gamma(t) $ is the dot product of two unit vectors: one in the
direction of $ {\bf r} $, and the other in the direction of $ \dot{\bf r} $.

From 
\begin{eqnarray}
|{\bf w}|^2 \hspace{-.2cm} & = & \hspace{-.2cm}
 |\dot{\bf r}|^2 - 2\left(\frac{\bf r}{|{\bf r}|} \cdot \dot{\bf r} \right)^2
 + \left(\frac{\bf r}{|{\bf r}|} \cdot \dot{\bf r} \right)^2 \nonumber \\
\hspace{-.2cm} & = & \hspace{-.2cm}
|\dot{\bf r}|^2 \left(1 - \Gamma^2 \right),
\end{eqnarray}
it follows that $ (1-\Gamma^2) $ is a measure of departure from motion 
camouflage.

\section{Feedback law for motion camouflage}

Using the planar setting as a guide, the curvature controls
to achieve motion camouflage in three dimensions can be
systematically derived.  Indeed, this is a major advantage of
representing trajectories using natural Frenet frames.  
However, for ease of exposition,
we instead begin by presenting the control law in an
intuitively appealing and biologically plausible form,
followed by the calculations demonstrating its effectiveness.

\subsection{Feedback law and interpretation}

Using the BAC-CAB identity, $ \tilde{\bf a} \times 
(\tilde{\bf b} \times \tilde{\bf c})
 = \tilde{\bf b}(\tilde{\bf a} \cdot \tilde{\bf c}) - 
 \tilde{\bf c} (\tilde{\bf a} \cdot \tilde{\bf b}) $, for
arbitrary vectors $ \tilde{\bf a} $, $ \tilde{\bf b} $, 
$ \tilde{\bf c} $, we observe that
\begin{equation}
\label{wcrossdefn}
{\bf w} = \dot{\bf r}\left(\frac{\bf r}{|{\bf r}|} \cdot
\frac{\bf r}{|{\bf r}|}\right) - \frac{\bf r}{|{\bf r}|} \left(
\frac{\bf r}{|{\bf r}|} \cdot \dot{\bf r} \right) 
 = \frac{\bf r}{|{\bf r}|} \times \left( \dot{\bf r} \times
\frac{\bf r}{|{\bf r}|} \right),
\end{equation}
\begin{eqnarray}
\label{wcrossr}
{\bf w} \times \frac{\bf r}{|{\bf r}|} \hspace{-.2cm} & = & \hspace{-.2cm}
  \left[\frac{\bf r}{|{\bf r}|}
 \times\left(\dot{\bf r} \times \frac{\bf r}{|{\bf r}|} \right)\right]
 \times \frac{\bf r}{|{\bf r}|} \nonumber \\
 \hspace{-.2cm} & = & \hspace{-.2cm}
- \frac{\bf r}{|{\bf r}|} \left[\left(\dot{\bf r} \times
 \frac{\bf r}{|{\bf r}|} \right) \cdot \frac{\bf r}{|{\bf r}|} \right]
 + \left(\dot{\bf r} \times \frac{\bf r}{|{\bf r}|} \right) \nonumber \\
\hspace{-.2cm} & = & \hspace{-.2cm}
  \dot{\bf r} \times \frac{\bf r}{|{\bf r}|},
\end{eqnarray}
and we conclude from (\ref{wcrossr}) that $ \left( \dot{\bf r} \times
{\bf r}/|{\bf r}| \right) $ is a biologically plausible quantity to
appear in a feedback law, since it only requires sensing 
$ {\bf w} $ and $ {\bf r}/|{\bf r}| $.

\begin{figure}
\hspace{.5cm}
\epsfxsize=6.7cm
\epsfbox{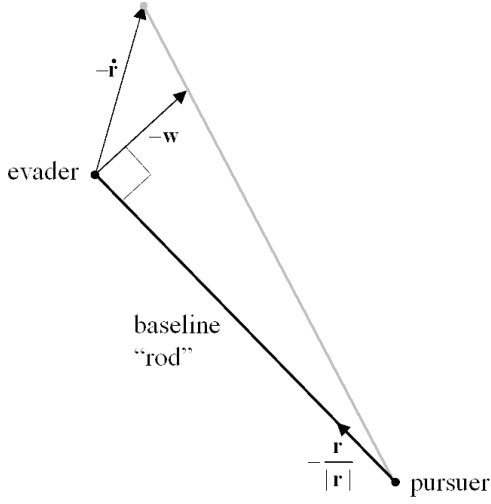}
\caption{\label{rodfig} Motion of the rod connecting the evader
to the pursuer, from the point of view of the pursuer.  The 
angular velocity of the rod, \boldmath$\omega$\unboldmath, is a
vector pointing into the page.}
\end{figure}

The quantity $ \left( \dot{\bf r} \times {\bf r}/|{\bf r}| \right) $ 
can be interpreted in terms of an angular-velocity-like quantity.  
From the point of view of the pursuer, consider an extensible
rod connecting the pursuer and evader positions.  The motion of
the evader (relative to the pursuer) contributes to change in the length
of this rod, as well as to angular velocity of the rod (viewed from
the pursuer - see figure \ref{rodfig}).  The transverse component
of the velocity of
the evader (viewed from the pursuer) is simply
\begin{eqnarray}
(\dot{\bf r}_e - \dot{\bf r}_p) - \left[(\dot{\bf r}_e - \dot{\bf r}_p)
 \cdot \frac{{\bf r}_e - {\bf r}_p}{|{\bf r}_e - {\bf r}_p|} \right]
 \frac{{\bf r}_e - {\bf r}_p}{|{\bf r}_e - {\bf r}_p|}
 \hspace{-6cm} & & \nonumber \\
 \hspace{-.2cm} & = & \hspace{-.2cm}
-\dot{\bf r} - \left[-\dot{\bf r} \cdot \left(-\frac{\bf r}{|{\bf r}|} 
 \right) \right]
 \left(-\frac{\bf r}{|{\bf r}|} \right) \nonumber \\
 \hspace{-.2cm} & = & \hspace{-.2cm}
 - {\bf w},
\end{eqnarray}
which can also be expressed as
\begin{equation}
\label{omegadefn}
-{\bf w} = \mbox{\boldmath$\omega$\unboldmath} \times (-{\bf r}),
\end{equation}
where \boldmath$\omega$ \unboldmath is the 
corresponding angular velocity of the rod.
From (\ref{wcrossdefn}) and (\ref{omegadefn}) we conclude that
\begin{equation}
\frac{\bf r}{|{\bf r}|} \times \left( \dot{\bf r} \times
\frac{\bf r}{|{\bf r}|} \right)=
\left(\frac{\bf r}{|{\bf r}|^2} \times \dot{\bf r} \right)
 \times {\bf r} =
  \mbox{\boldmath$\omega$\unboldmath} \times {\bf r} \nonumber \\
\end{equation}
and hence
\begin{equation}
\label{omegaeqn}
\mbox{\boldmath$\omega$\unboldmath} =
 \frac{\bf r}{|{\bf r}|^2} \times \dot{\bf r}.
\end{equation}
Thus, the quantity  $ \left( \dot{\bf r} \times {\bf r}/|{\bf r}| \right) $
is simply $-\mbox{\boldmath$\omega$\unboldmath}$ scaled by $ |{\bf r}| $.

For convenience in the calculations below, we define
\begin{equation}
\label{adefn}
{\bf a} = {\bf x}_p \times  
 \left(\dot{\bf r} \times \frac{\bf r}{|{\bf r}|} \right),
\end{equation}
and express the feedback law as
\begin{eqnarray}
\label{updefn}
u_p \hspace{-.2cm} & = & \hspace{-.2cm}
 \mu ({\bf a} \cdot {\bf y}_p), \\
\label{vpdefn}
v_p \hspace{-.2cm} & = & \hspace{-.2cm}
 \mu ({\bf a} \cdot {\bf z}_p),
\end{eqnarray}
where $ \mu > 0 $ is a constant feedback gain. 
The quantity $ \mu \nu_p^2 {\bf a} $ can then be interpreted as the
lateral component of the acceleration vector of the pursuer.  
Consistent with the
fact that $ u_p $ and $ v_p $ can only change the direction
of the pursuer's motion and not its speed, we note that
$ \mu \nu_p^2 {\bf a} $ is transverse to the
direction of motion of the pursuer, $ {\bf x}_p $: i.e.,
$ {\bf a} \cdot {\bf x}_p = 0 $.
 
 Using the formula
$ \tilde{\bf a} \cdot (\tilde{\bf b} \times \tilde{\bf c}) = 
\tilde{\bf b} \cdot (\tilde{\bf c} \times
\tilde{\bf a}) $ for the scalar triple product, 
where $ \tilde{\bf a} $, $ \tilde{\bf b} $, $ \tilde{\bf c} $ are
arbitrary vectors, we compute
\begin{eqnarray}
\label{upformula}
u_p \hspace{-.2cm} & = & \hspace{-.2cm}
 \mu \left[ {\bf x}_p \times
 \left(\dot{\bf r} \times \frac{\bf r}{|{\bf r}|} \right) \right] 
 \cdot {\bf y}_p \nonumber \\
 \hspace{-.2cm} & = & \hspace{-.2cm}
 \mu\left[ \left(\dot{\bf r} \times \frac{\bf r}{|{\bf r}|} \right)
 \cdot ({\bf y}_p \times {\bf x}_p) \right] \nonumber \\
 \hspace{-.2cm} & = & \hspace{-.2cm}
 -\mu \left[ \left(\dot{\bf r} \times \frac{\bf r}{|{\bf r}|} \right)
 \cdot {\bf z}_p \right],
\end{eqnarray}
and similarly,
\begin{equation}
v_p = \mu\left[\left(\dot{\bf r} \times \frac{\bf r}{|{\bf r}|}
 \right) \cdot {\bf y}_p \right].
\end{equation}

\vspace{.5cm}

\noindent
{\bf Remark}: It is easy to see that in the planar setting, we
recover the planar steering law for motion camouflage presented
in \cite{mcarxiv05}.  If $ {\bf x}_p $, $ {\bf x}_e $, and $ {\bf r} $
all lie in the same plane, then $ \dot{\bf r} $ also lies in that plane,
and (\ref{upformula}) becomes
\begin{equation}
u_p =
 - \mu\left[\left(\dot{\bf r} \times \frac{\bf r}{|{\bf r}|} \right) \cdot
 {\bf z}_p \right] =
 - \mu\left(\frac{\bf r}{|{\bf r}|} \cdot \dot{\bf r}^{\perp} \right),
\end{equation}
where the notion $ {\bf q}^{\perp} $ represents the vector $ {\bf q} $
rotated counterclockwise in the plane by $ \pi/2 $.  Furthermore,
without loss of generality, we identify
$ {\bf y}_p $ with $ {\bf x}_p^\perp $, and
$ {\bf z}_p $ with the unit vector perpendicular to
the plane of motion. $ \Box $

\subsection{Behavior of $ \Gamma $ under the feedback law}

Differentiating $ \Gamma $ along trajectories of (\ref{pursuer3d})
and (\ref{evader3d}) gives
\begin{eqnarray}
\label{dotgamma}
\dot{\Gamma} \hspace{-.2cm} & = & \hspace{-.3cm}
\left(\frac{\dot{\bf r} \cdot \dot{\bf r} + {\bf r} \cdot \ddot{\bf r}}
{|{\bf r}||\dot{\bf r}|} \right)
 - \left(\frac{{\bf r}\cdot \dot{\bf r}}{|\dot{\bf r}|} \right) \hspace{-.1cm}
 \left( \frac{{\bf r} \cdot \dot{\bf r}}{|{\bf r}|^3} \right)
 - \left( \frac{{\bf r}\cdot \dot{\bf r}}{|{\bf r}|} \right) \hspace{-.1cm}
 \left( \frac{\dot{\bf r} \cdot \ddot{\bf r}}{|\dot{\bf r}|^3} \right)
\nonumber \\
\hspace{-.2cm} & = & \hspace{-.3cm}
\frac{|\dot{\bf r}|}{|{\bf r}|} \hspace{-.1cm}
 \left[1 - {\left(\frac{\bf r}{|{\bf r}|}
 \cdot \frac{\dot{\bf r}}{|\dot{\bf r}|}\right) \hspace{-.1cm}}^2\right]
  \hspace{-.05cm} + \hspace{-.05cm} 
 \frac{1}{|\dot{\bf r}|}\left[ \frac{\bf r}{|{\bf r}|} -
 \left(\frac{\bf r}{|{\bf r}|}\cdot \frac{\dot{\bf r}}{|\dot{\bf r}|}\right)
 \frac{\dot{\bf r}}{|\dot{\bf r}|} \right]
 \hspace{-.1cm} \cdot \hspace{-.05cm} \ddot{\bf r}. \nonumber \\
\end{eqnarray}
We also have
\begin{eqnarray}
\ddot{\bf r} \hspace{-.2cm} & = & \hspace{-.2cm}
 \dot{\nu}_p {\bf x}_p - \dot{\nu}_e {\bf x}_e
 + \nu_p \dot{\bf x}_p - \nu_e \dot{\bf x}_e \nonumber \\
\hspace{-.2cm} & = & \hspace{-.2cm}
 \dot{\nu}_p {\bf x}_p - \dot{\nu_e} {\bf x}_e
 + \nu_p^2 ({\bf y}_p u_p + {\bf z}_p v_p) 
 - \nu_e^2 ({\bf y}_e u_e + {\bf z}_e v_e). \nonumber \\
\end{eqnarray} 
If we define
\begin{equation}
{\bf b} = \frac{1}{|\dot{\bf r}|} \left[ \frac{\bf r}{|{\bf r}|} -
 \left(\frac{\bf r}{|{\bf r}|}\cdot \frac{\dot{\bf r}}{|\dot{\bf r}|}\right)
 \frac{\dot{\bf r}}{|\dot{\bf r}|} \right],
\end{equation}
then
\begin{eqnarray}
{\bf b} \cdot \ddot{\bf r} \hspace{-.2cm} & = & \hspace{-.2cm}
 \dot{\nu}_p ({\bf b} \cdot {\bf x}_p) - \dot{\nu}_e({\bf b} \cdot {\bf x}_e)
 \nonumber \\ & &
 + \nu_p^2 \left[({\bf b} \cdot {\bf y}_p) u_p +
 ({\bf b} \cdot {\bf z}_p) v_p \right] \nonumber \\ & &
 - \nu_e^2 \left[({\bf b} \cdot {\bf y}_e) u_e + 
 ({\bf b} \cdot {\bf z}_e) v_e \right],
\end{eqnarray}
and the only term of $ \dot{\Gamma} $ into which the controls $ u_p $ 
and $ v_p $ explicitly enter is
\begin{equation}
\nu_p^2 \left[({\bf b} \cdot {\bf y}_p) u_p +
 ({\bf b} \cdot {\bf z}_p) v_p \right].
\end{equation}
Using (\ref{updefn}) and (\ref{vpdefn}),
\begin{eqnarray}
\nu_p^2 \left[({\bf b} \cdot {\bf y}_p) u_p +
 ({\bf b} \cdot {\bf z}_p) v_p \right] \hspace{-3cm} & & \nonumber \\
\hspace{-.2cm} & = & \hspace{-.2cm}  \mu \;
\nu_p^2 \left[({\bf b} \cdot {\bf y}_p) \left({\bf a} \cdot {\bf y}_p \right) +
 ({\bf b} \cdot {\bf z}_p) \left({\bf a} \cdot {\bf z}_p \right) \right]
 \nonumber \\
\hspace{-.2cm} & = & \hspace{-.2cm}
  \mu \;
\nu_p^2 \left[({\bf b} \cdot {\bf a})-({\bf b} \cdot {\bf x}_p)
 ({\bf a} \cdot {\bf x}_p) \right]  \nonumber \\
\hspace{-.2cm} & = & \hspace{-.2cm}
  \mu \;
\nu_p^2 ({\bf b} \cdot {\bf a}),
\end{eqnarray}
where we have also used the identity
\begin{equation}
{\bf a} \cdot {\bf b} =
({\bf a} \cdot {\bf x}_p)({\bf b} \cdot {\bf x}_p)
 +({\bf a} \cdot {\bf y}_p)({\bf b} \cdot {\bf y}_p)
 +({\bf a} \cdot {\bf z}_p)({\bf b} \cdot {\bf z}_p),
\end{equation}
and $ {\bf a} \cdot {\bf x}_p = 0 $.

Using the BAC-CAB identity, we observe that
\begin{equation}
{\bf b} = \frac{1}{|\dot{\bf r}|} \hspace{-.1cm}
 \left[\frac{\dot{\bf r}}{|\dot{\bf r}|} \hspace{-.05cm} \times \hspace{-.05cm}
 \left(\frac{\bf r}{|{\bf r}|} \hspace{-.05cm} \times \hspace{-.05cm}
 \frac{\dot{\bf r}}{|\dot{\bf r}|}
 \right) \hspace{-.05cm} \right] 
 = - \frac{1}{|\dot{\bf r}|^3} \hspace{-.05cm}
  \left[\dot{\bf r}  \times \hspace{-.05cm}
 \left(\dot{\bf r}  \times \hspace{-.05cm}
  \frac{\bf r}{|{\bf r}|} \right) \hspace{-.05cm} \right],
\end{equation}
so that
\begin{equation}
{\bf b} \cdot {\bf a} =  - \frac{1}{|\dot{\bf r}|^3} \left[\dot{\bf r} \times
 \left(\dot{\bf r} \times   \frac{\bf r}{|{\bf r}|} \right) \right]
 \cdot \left[{\bf x}_p \times
 \left(\dot{\bf r} \times \frac{\bf r}{|{\bf r}|} \right) \right].
\end{equation}
Using the identity $ (\tilde{\bf a} \times \tilde{\bf b}) \cdot 
(\tilde{\bf c} \times \tilde{\bf d})=
(\tilde{\bf a} \cdot \tilde{\bf c})(\tilde{\bf b}\cdot \tilde{\bf d})
-(\tilde{\bf a}\cdot \tilde{\bf d})
 (\tilde{\bf b} \cdot \tilde{\bf c}) $, for arbitrary vectors
$ \tilde{\bf a} $, $ \tilde{\bf b} $, $ \tilde{\bf c} $, $ \tilde{\bf d} $, and
\begin{equation}
\left|\dot{\bf r} \times \frac{\bf r}{|{\bf r}|}\right|^2 =
 |\dot{\bf r}|^2\left(1-\Gamma^2\right),
\end{equation}
we compute
\begin{eqnarray}
{\bf b} \cdot {\bf a} \hspace{-.2cm} & = & \hspace{-.2cm}
- \frac{1}{|\dot{\bf r}|} \left(\dot{\bf r} \cdot {\bf x}_p\right)
 \left(1-\Gamma^2\right) \nonumber \\ & & 
 +\frac{1}{|\dot{\bf r}|^3} \left[\dot{\bf r} \cdot  
 \left(\dot{\bf r} \times \frac{\bf r}{|{\bf r}|} \right)\right]
 \left[ \left(\dot{\bf r} \times \frac{\bf r}{|{\bf r}|} \right) \cdot 
 {\bf x}_p \right] \nonumber \\
\hspace{-.2cm} & = & \hspace{-.2cm}
- \left(1-\Gamma^2\right)
 \left(\frac{\dot{\bf r}}{|\dot{\bf r}|} \cdot {\bf x}_p \right).
\end{eqnarray}

\vspace{.5cm}

\noindent
{\bf Remark}: For the foregoing calculations to make sense, we
require $ |{\bf r}|>0 $ and $ |\dot{\bf r}|>0 $.  The condition
$ |{\bf r}|> 0 $ is a non-collision condition,
and does not pose any difficulty for us
because our analysis of approach to the state of motion camouflage
takes place away from the collision state.  Later, we will
impose hypotheses that also ensure  $ |\dot{\bf r}| > 0 $
for all time.  (Note that in the constant-speed setting, 
$ 0 < \nu_e/\nu_p < 1 $ is sufficient to ensure $ |\dot{\bf r}|> 0 $
\cite{mcarxiv05}.) $ \Box $

\vspace{.5cm}

\noindent
{\bf Remark}: Provided $ \nu_p > \nu_e $, we have
\begin{equation}
\left(\frac{\dot{\bf r}}{|\dot{\bf r}|} \cdot {\bf x}_p \right)
 = \frac{1}{|\dot{\bf r}|}\left[\nu_p - \nu_e\left({\bf x}_p \cdot {\bf x}_e
 \right) \right] > 0,
\end{equation}
so that
\begin{equation}
{\bf b} \cdot {\bf a} \le 0,
\end{equation}
and therefore the only term in $ \dot{\Gamma} $ explicitly 
involving the controls $ u_p $ and $ v_p $ satisfies
\begin{equation}
\nu_p^2 \left[({\bf b} \cdot {\bf y}_p) u_p +
 ({\bf b} \cdot {\bf z}_p) v_p \right] \le 0. 
\end{equation}
$ \Box $

\vspace{.5cm}

To summarize, (\ref{dotgamma}) becomes
\begin{eqnarray}
\label{dotgamma2}
\dot{\Gamma} \hspace{-.2cm} & = & \hspace{-.2cm}
- \left(1-\Gamma^2\right) \left[
 \frac{\mu \nu_p^2}{|\dot{\bf r}|}
 \left(\nu_p - \nu_e\left({\bf x}_p \cdot {\bf x}_e \right) \right)
 - \frac{|\dot{\bf r}|}{|{\bf r}|} \right]\nonumber \\ & &
 +  \dot{\nu}_p ({\bf b} \cdot {\bf x}_p) - \dot{\nu}_e({\bf b} \cdot {\bf x}_e)
 \nonumber \\ & &
 - \nu_e^2 \left[({\bf b} \cdot {\bf y}_e) u_e +
 ({\bf b} \cdot {\bf z}_e) v_e \right].
\end{eqnarray}
Noting that 
\begin{equation}
|{\bf b}|^2 = 
\frac{1}{|\dot{\bf r}|^2} \left(1 - \Gamma^2 \right),
\end{equation}
we see that
\begin{eqnarray}
\big|\nu_e^2 \left[({\bf b} \cdot {\bf y}_e) u_e +
 ({\bf b} \cdot {\bf z}_e) v_e \right] \big| \hspace{-4cm} & & \nonumber \\
\hspace{-.2cm} & \le & \hspace{-.2cm}
\frac{\nu_e^2}{|\dot{\bf r}|} \sqrt{1 - \Gamma^2} 
 \max \left(\sqrt{u_e^2+v_e^2}\right),
\end{eqnarray}
where $  \max \left(\sqrt{u_e^2+v_e^2}\right) $ is an a priori bound
on the maximum absolute curvature of the evader trajectory.
Similarly,
\begin{eqnarray}
\big| \dot{\nu}_p ({\bf b} \cdot {\bf x}_p) 
 - \dot{\nu}_e({\bf b} \cdot {\bf x}_e) \big| \hspace{-3cm}
& & \nonumber \\
\hspace{-.2cm} & \le & \hspace{-.2cm}
\frac{1}{|\dot{\bf r}|} \sqrt{1 - \Gamma^2} 
 \left(|\dot{\nu}_p| + | \dot{\nu}_e | \right) 
\nonumber \\
\hspace{-.2cm} & \le & \hspace{-.2cm}
\frac{1}{|\dot{\bf r}|} \sqrt{1 - \Gamma^2} 
 \left(\alpha_p + \alpha_e\right),
\end{eqnarray}
where $ \alpha_p $ is an upper bound on $ |\dot{\nu}_p| $,
and $ \alpha_e $ is an upper bound on $ |\dot{\nu}_e| $.
From (\ref{dotgamma2}) we then conclude
\begin{eqnarray}
\label{dotgamma3}
\dot{\Gamma} \hspace{-.2cm} & \le & \hspace{-.2cm}
- \left(1-\Gamma^2\right) \left[
 \frac{\mu \nu_p^2}{|\dot{\bf r}|}
 \left(\nu_p - \nu_e\left({\bf x}_p \cdot {\bf x}_e \right) \right)
 - \frac{|\dot{\bf r}|}{|{\bf r}|} \right]\nonumber \\ & &
+ \frac{1}{|\dot{\bf r}|} \sqrt{1 - \Gamma^2}
 \left[\alpha_p + \alpha_e + \nu_e^2\max \left(\sqrt{u_e^2+v_e^2}\right)
 \right].  \nonumber \\ 
\end{eqnarray}

\subsection{Bounds and estimates}

Having bounded $ \dot{\Gamma} $ as in (\ref{dotgamma3}),
we proceed in analogy with the planar setting \cite{mcarxiv05}.
We hypothesize that a constant
$ \nu_{\mathit max} $ exists such that
\begin{equation}
 \frac{\nu_e}{\nu_p} \le \nu_{\mathit max} < 1,
\end{equation}
for all time.  We also assume that constants
$ \nu_p^{\mathit low} $, $ \nu_p^{\mathit high} $,
$ \nu_e^{\mathit low} $, and $ \nu_e^{\mathit high} $ exist such that
\begin{eqnarray}
0 \hspace{-.2cm} & < & \hspace{-.2cm} 
 \nu_p^{\mathit low} \le \nu_p \le \nu_p^{\mathit high} < \infty, \\
0 \hspace{-.2cm} & < & \hspace{-.2cm}
  \nu_e^{\mathit low} \le \nu_e \le \nu_e^{\mathit high} < \infty,
\end{eqnarray}
for all time, and observe that 
\begin{equation}
\label{dotrbound}
0 < \nu_p^{\mathit low}(1-\nu_{\mathit max}) \le |\dot{\bf r}|
  \le \nu_p^{\mathit high}
 (1+\nu_{\mathit max}).
\end{equation}

We define the constant $ c_1 > 0 $ as
\begin{equation}
c_1 = \frac{\left[\alpha_p + \alpha_e + 
 (\nu_e^{\mathit high})^2\max \left(\sqrt{u_e^2+v_e^2}\right)\right]}
{\nu_p^{\mathit low}(1-\nu_{\mathit max})}.
\end{equation}
Given $ \mu > 0 $ sufficiently large and $ r_o > 0 $, we define $ c_0 > 0 $ by
\begin{equation}
c_0 = \left(\frac{(\nu_p^{\mathit low})^3
 (1-\nu_{\mathit max})}{\nu_p^{\mathit high}
 (1+\nu_{\mathit max})} \right)\mu 
 - \frac{\nu_p^{\mathit high}(1+\nu_{\mathit max})}{r_o},
\end{equation} 
so that
\begin{equation}
\label{mudecomp}
\mu  = \left(\frac{\nu_p^{\mathit high}
 (1+\nu_{\mathit max})}{(\nu_p^{\mathit low})^3(1-\nu_{\mathit max})}
 \right) \left(
 \frac{\nu_p^{\mathit high}(1+\nu_{\mathit max})}{r_o} + c_0 \right),
\end{equation}
and hence
\begin{equation}
\mu \ge \hspace{-.1cm} \left(\frac{\nu_p^{\mathit high}
 (1+\nu_{\mathit max})}{(\nu_p^{\mathit low})^3(1-\nu_{\mathit max})}
 \right) \hspace{-.15cm} \left(
 \frac{\nu_p^{\mathit high}(1+\nu_{\mathit max})}
 {|{\bf r}|} + c_0 \hspace{-.05cm} \right), 
\end{equation}
$ \forall |{\bf r}| \ge r_o. $ 
Thus, for $ |{\bf r}| \ge r_o $, (\ref{dotgamma3}) becomes
\begin{eqnarray}
\dot{\Gamma} \hspace{-.3cm} & \le & \hspace{-.3cm}
   - \hspace{-.05cm}
 \left(1 \hspace{-.05cm} -  \hspace{-.05cm} \Gamma^2 \right) \hspace{-.15cm}
 \Bigg[ \hspace{-.15cm} \left( \hspace{-.1cm} \frac{\nu_p^{\mathit high}
 (1 \hspace{-.05cm} + \hspace{-.05cm} \nu_{\mathit max})}
 {(\nu_p^{\mathit low})^3
 (1  \hspace{-.05cm} - \hspace{-.05cm} \nu_{\mathit max})}
 \hspace{-.1cm} \right)
 \hspace{-.15cm} \left( \hspace{-.1cm}
  \frac{\nu_p^{\mathit high}(1 \hspace{-.05cm} + \hspace{-.05cm}
 \nu_{\mathit max})}{|{\bf r}|}
  \hspace{-.05cm} +  \hspace{-.05cm} c_0 \hspace{-.1cm} \right)
 \nonumber \\ & & \hspace{1.3cm} \times \hspace{-.1cm}
 \left( \hspace{-.1cm}
 \frac{(\nu_p^{\mathit low})^3
 (1  \hspace{-.05cm} - \hspace{-.05cm} \nu_{\mathit max})}{\nu_p^{\mathit high}
 (1  \hspace{-.05cm} + \hspace{-.05cm} \nu_{\mathit max})} 
 \hspace{-.1cm} \right) \hspace{-.05cm}
 - \frac{\nu_p^{\mathit high}(1 \hspace{-.05cm} + \hspace{-.05cm}
 \nu_{\mathit max})}{|{\bf r}|} 
 \hspace{-.05cm} \Bigg] \nonumber \\ & &
 + \left(\sqrt{1-\Gamma^2}\right)c_1
  \nonumber \\
 \hspace{-.3cm} & = & \hspace{-.3cm}
   - \left(1-\Gamma^2 \right) c_0 + \left(\sqrt{1-\Gamma^2} \right)c_1.
\end{eqnarray}
Suppose that given $ 0<\epsilon << 1$, we take $ \mu > 0 $ sufficiently
large so that there exists $ c_0 $ satisfying
$ c_0 \ge 2c_1/\sqrt{\epsilon} $.  Then for $ (1-\Gamma^2) > \epsilon $,
\begin{eqnarray}
\label{dotgammabound2}
\dot{\Gamma}  \hspace{-.2cm} & \le & \hspace{-.25cm} 
 - \left(1-\Gamma^2 \right) c_0 + \left(\sqrt{1-\Gamma^2} \right)c_1
 \nonumber \\
 \hspace{-.2cm} & = & \hspace{-.2cm}
  - \left(1-\Gamma^2 \right) \left( c_0
 - \frac{c_1}{\sqrt{1-\Gamma^2}} \right) \nonumber \\
  \hspace{-.2cm} & \le & \hspace{-.25cm}
  - \left(1-\Gamma^2 \right) \left( c_0
 - \frac{c_1}{\sqrt{\epsilon}}\right) \nonumber \\
 \hspace{-.2cm} & = & \hspace{-.2cm}
 -  \left(1-\Gamma^2 \right) c_2,
\end{eqnarray} 
where
\begin{equation}
\label{c2defn}
c_2 = c_0 - \frac{c_1}{\sqrt{\epsilon}} > 0.
\end{equation}

\vspace{.5cm}

\noindent
{\bf Remark}: There are two possibilities for 
\begin{equation}
\label{oneminusgammasqreps}
 (1-\Gamma^2) \le \epsilon.
\end{equation} 
The state we seek to drive the system toward has $ \Gamma \approx -1 $;
however, (\ref{oneminusgammasqreps}) can also be satisfied for 
$ \Gamma \approx 1 $.  (Recall that $ -1 \le \Gamma \le 1 $.)
There is always a set of initial conditions
such that (\ref{oneminusgammasqreps}) is satisfied with $ \Gamma \approx 1 $.
We can address this issue as follows: let $ \epsilon_o > 0 $ denote how
close to $ -1 $ we wish to drive $ \Gamma $, and let $ \Gamma_0 = \Gamma(0) $
denote the initial value of $ \Gamma $.  Take
\begin{equation}
\epsilon = \min(\epsilon_o,1-\Gamma_0^2),
\end{equation}
so that (\ref{dotgammabound2}) with (\ref{c2defn})
applies from time $ t = 0 $. $ \Box $

\vspace{.5cm}

From (\ref{dotgammabound2}), we can write
\begin{equation}
\frac{d\Gamma}{1-\Gamma^2} \le - c_2 dt,
\end{equation}
which, integrating both sides, leads to
\begin{equation}
\int_{\Gamma_0}^{\Gamma} \frac{d\tilde{\Gamma}}{1-\tilde{\Gamma}^2} \le 
 -c_2 \int_0^t d\tilde{t} = -c_2 t,
\end{equation}
where $ \Gamma_0 = \Gamma(t=0) $.
Noting that
\begin{equation}
\int_{\Gamma_0}^{\Gamma}  \frac{d\tilde{\Gamma}}{1-\tilde{\Gamma}^2}
 = \int_{\Gamma_0}^{\Gamma} d(\tanh^{-1} \tilde{\Gamma}) 
 = \tanh^{-1} \Gamma - \tanh^{-1} \Gamma_0,
\end{equation}
we see that for $ |{\bf r}| \ge r_o $, (\ref{dotgammabound2}) implies
\begin{equation}
\label{gammabound}
\Gamma(t) \le \tanh\left(\tanh^{-1}\Gamma_0 - c_2 t \right),
\end{equation}
where we have used the fact that $ \tanh^{-1}(\cdot) $ is
a monotone increasing function.

Now we consider estimating how long $ |{\bf r}| \ge r_o $, 
which in turn determines how large $ t $ can become in inequality
(\ref{gammabound}), and hence how close to $-1$ will $ \Gamma(t) $ be 
driven.  From (\ref{ddtnormr}) and (\ref{gammadotprod}) we have
\begin{equation}
\frac{d}{dt}|{\bf r}| = \Gamma(t) |\dot{\bf r}|,
\end{equation}
which from (\ref{dotrbound}) and $ |\Gamma(t)| \le 1 $, $ \forall t $,
implies 
\begin{equation}
\label{dotnormrbound}
\frac{d}{dt}|{\bf r}| \ge -|\Gamma(t)|\nu_p^{\mathit high}(1+\nu_{\mathit max}) 
 \ge - \nu_p^{\mathit high}(1+\nu_{\mathit max}).
\end{equation}
From (\ref{dotnormrbound}), we conclude that
\begin{equation}
|{\bf r}(t)| \ge |{\bf r}(0)| - \nu_p^{\mathit high}(1+\nu_{\mathit max})t, 
 \;\; \forall t \ge 0,
\end{equation}
and, more to the point,
\begin{equation}
\label{dotnormrbound2}
|{\bf r}(t)| \ge r_o, \;\; \forall t \le \frac{|{\bf r}(0)|-r_o}
 {\nu_p^{\mathit high}(1+\nu_{\mathit max})}.
\end{equation}
For (\ref{dotnormrbound2}) to be meaningful for the problem at hand,
we assume that $ |{\bf r}(0)| > r_o $.  Then defining
\begin{equation}
\label{bigtdefn}
T = \frac{|{\bf r}(0)|-r_o}{\nu_p^{\mathit high}(1+\nu_{\mathit max})} > 0 
\end{equation}
to be the minimum interval of time over which we can guarantee
that $ \dot{\Gamma} \le 0 $, we conclude that
\begin{equation}
\label{gammafinal}
\Gamma(T)  \le  
  \tanh \left( \tanh^{-1} \Gamma_0 - c_2 T \right). 
\end{equation}

From (\ref{gammafinal}), we see that by choosing $ c_2 $ sufficiently
large (which can be accomplished by choosing $ c_0 \ge 2c_1/\sqrt{\epsilon} $
sufficiently large), we can force $ \Gamma(T) \le -1+\epsilon $.
Noting that
\begin{equation}
\tanh(x) 
 \le -1+\epsilon
 \Longleftrightarrow x \le \frac{1}{2} \ln 
 \left(\frac{\epsilon}{2 -\epsilon} \right),
\end{equation}
for $ 0 < \epsilon << 1 $,
we see that 
\begin{equation}
\Gamma(T) \le -1 + \epsilon
  \Longleftrightarrow
\tanh^{-1} \Gamma_0 - c_2 T \le  \frac{1}{2} \ln 
 \left(\frac{\epsilon}{2 -\epsilon} \right).
\end{equation}
Thus, if $ c_0 \ge 2c_1/\sqrt{\epsilon} $ is taken to be sufficiently
large that
\begin{equation}
\label{c2bound}
c_2 \ge \nu_p^{\mathit high}(1+\nu_{\mathit max})
 \frac{\tanh^{-1} \Gamma_0 - \frac{1}{2} \ln
 \left(\frac{\epsilon}{2 -\epsilon} \right)}
 { |{\bf r}(0)| - r_o },
\end{equation}
then we are guaranteed (under the conditions mentioned in the 
above calculations) to achieve $ \Gamma(t_1) \le -1 + \epsilon $
at some finite time $ t_1 \le T $.

\subsection{Statement of result}

\noindent
{\bf Definition} \cite{mcarxiv05}: Given the system (\ref{pursuer3d}),
(\ref{evader3d}) with $ \Gamma $ defined by (\ref{gammadefn}), we say
that ``motion camouflage is accessible in finite time'' if for any
$ \epsilon > 0 $ there exists a time $ t_1 > 0 $ such that
$ \Gamma(t_1) \le -1 + \epsilon $. 

\vspace{.5cm}

\noindent
{\bf Proposition}: Consider the system (\ref{pursuer3d}),
(\ref{evader3d}) with $ \Gamma $ defined by (\ref{gammadefn}) and
control law given by (\ref{adefn}) - (\ref{vpdefn}), with
the following hypotheses:
\begin{itemize}
\item[$ \hspace{-.05cm} (\hspace{-.05cm} A1) \hspace{-.05cm} $] 
 $ 0 < \nu_p^{\mathit low} \le
 \nu_p \le \nu_p^{\mathit high} < \infty $, where
  $ \nu_p^{\mathit low} $ and $ \nu_p^{\mathit high} $ are constants,
\item[$ \hspace{-.05cm} (\hspace{-.05cm} A2) \hspace{-.05cm} $] 
 $ 0 < \nu_e^{\mathit low} \le
 \nu_e \le \nu_e^{\mathit high} < \infty $, where
  $ \nu_e^{\mathit low} $ and $ \nu_e^{\mathit high} $ are constants,
\item[$ \hspace{-.05cm} (\hspace{-.05cm} A3) \hspace{-.05cm} $] $ 
 \nu_e/\nu_p \le \nu_{\mathit max} < 1 $, where
  $ \nu_{\mathit max} $ is constant,
\item[$ \hspace{-.05cm} (\hspace{-.05cm} A4) \hspace{-.05cm} $] 
 $ u_e $ and $ v_e $ are piecewise continuous and $ \sqrt{u_e^2+v_e^2} $ 
is bounded,
\item[$ \hspace{-.05cm} (\hspace{-.05cm} A5) \hspace{-.05cm} $] 
 $ \dot{\nu}_e $ and $ \dot{\nu}_p $ are piecewise continuous, 
 $ |\dot{\nu}_p| < \alpha_p $,
and $ |\dot{\nu}_e| < \alpha_e $, where $ \alpha_p $ and $ \alpha_e $ are
finite constants,
\item[$ \hspace{-.05cm} (\hspace{-.05cm} A6) \hspace{-.05cm} $] 
 $ \Gamma_0 = \Gamma(0) < 1 $, and
\item[$ \hspace{-.05cm} (\hspace{-.05cm} A7) \hspace{-.05cm} $] 
$ |{\bf r}(0)| > 0 $.
\end{itemize}
Motion camouflage is accessible in finite time
using high-gain feedback (i.e., by choosing $ \mu > 0 $ sufficiently
large).

\vspace{.5cm}

\noindent
{\bf Proof}: Analogous to the corresponding proof in \cite{mcarxiv05}.
Choose $ r_o > 0 $ such that $ r_o < |{\bf r}(0)| $.
Choose $ c_2 > 0 $ sufficiently large so as to satisfy (\ref{c2bound}),
and choose $ c_0 $ accordingly to ensure that (\ref{dotgammabound2}) holds
for $ \Gamma > -1 + \epsilon $.  Then defining $ \mu $ according to
(\ref{mudecomp}) ensures that $ \Gamma(T) \le -1 + \epsilon $,
where $ T > 0 $ is defined by (\ref{bigtdefn}).  $ \Box $

\section{Simulation Results}

\begin{figure}
\epsfxsize=8.5cm
\epsfbox{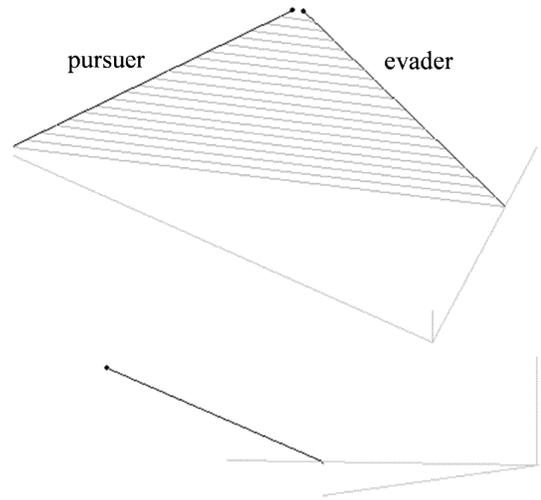}
\caption{\label{simfig1} Straight-line evader trajectory, and corresponding
pursuer trajectory. The pursuer and evader trajectories are the dark
lines (with dots at the final positions when the simulation is stopped).
The light lines connecting the pursuer and evader trajectories are 
baselines drawn at equally spaced time intervals.  The upper plot is
the view perpendicular to the baseline direction, and the lower plot is the
view along the baseline direction (so that the pursuer and evader
trajectories overlap). } 
\end{figure}

\begin{figure}
\epsfxsize=8.5cm
\epsfbox{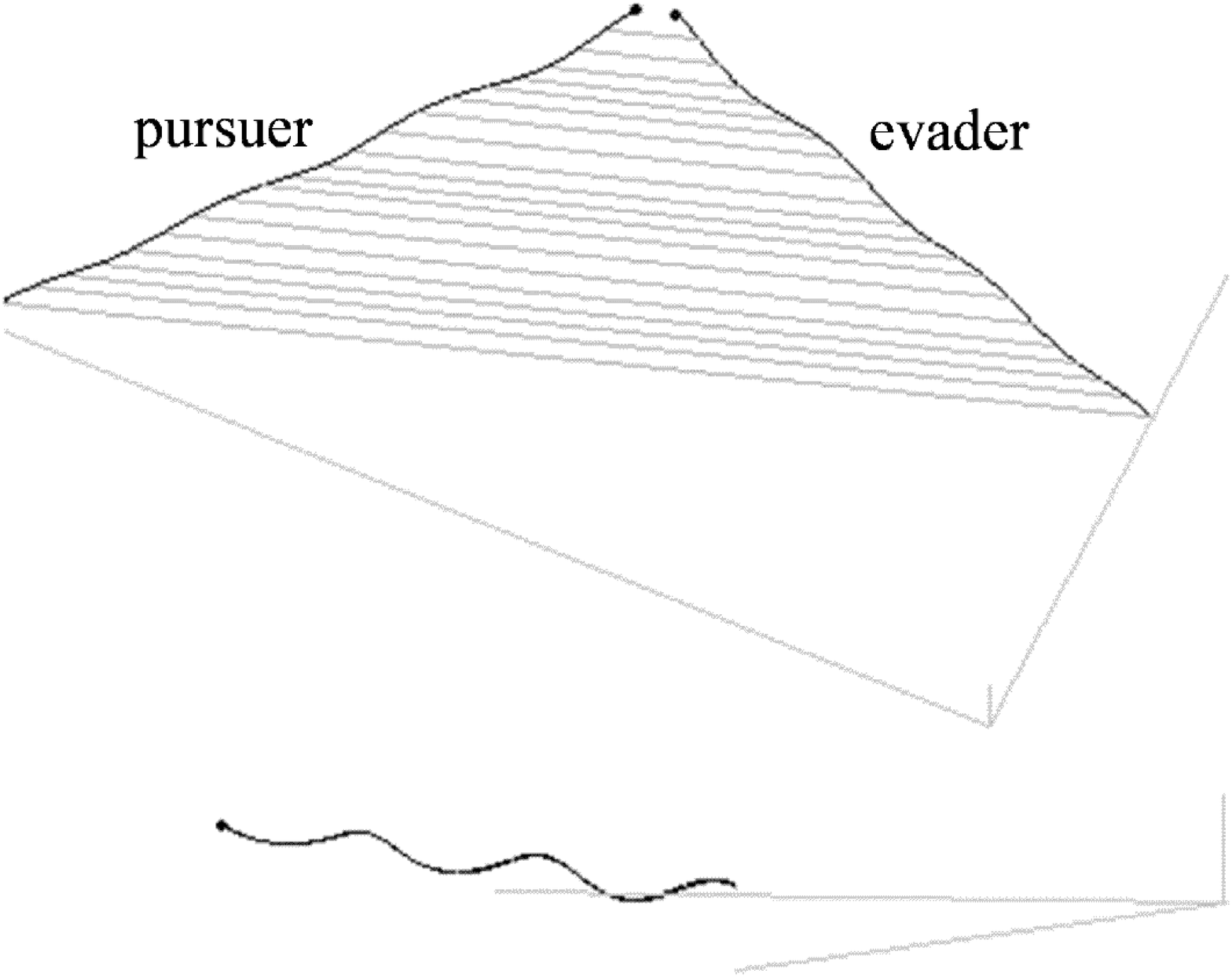}
\caption{\label{simfig2} Evader trajectory with sinusoidally
varying curvature inputs, and corresponding
pursuer trajectory.}  
\end{figure}

\begin{figure}
\epsfxsize=8.5cm
\epsfbox{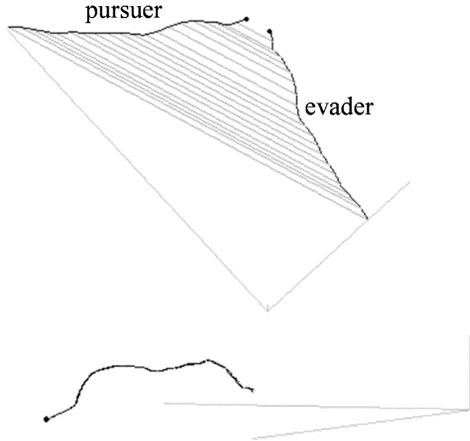}
\caption{\label{simfig3} Evader trajectory with randomly varying curvature
inputs, and corresponding
pursuer trajectory.}
\end{figure}

\begin{figure}
\epsfxsize=8.5cm
\epsfbox{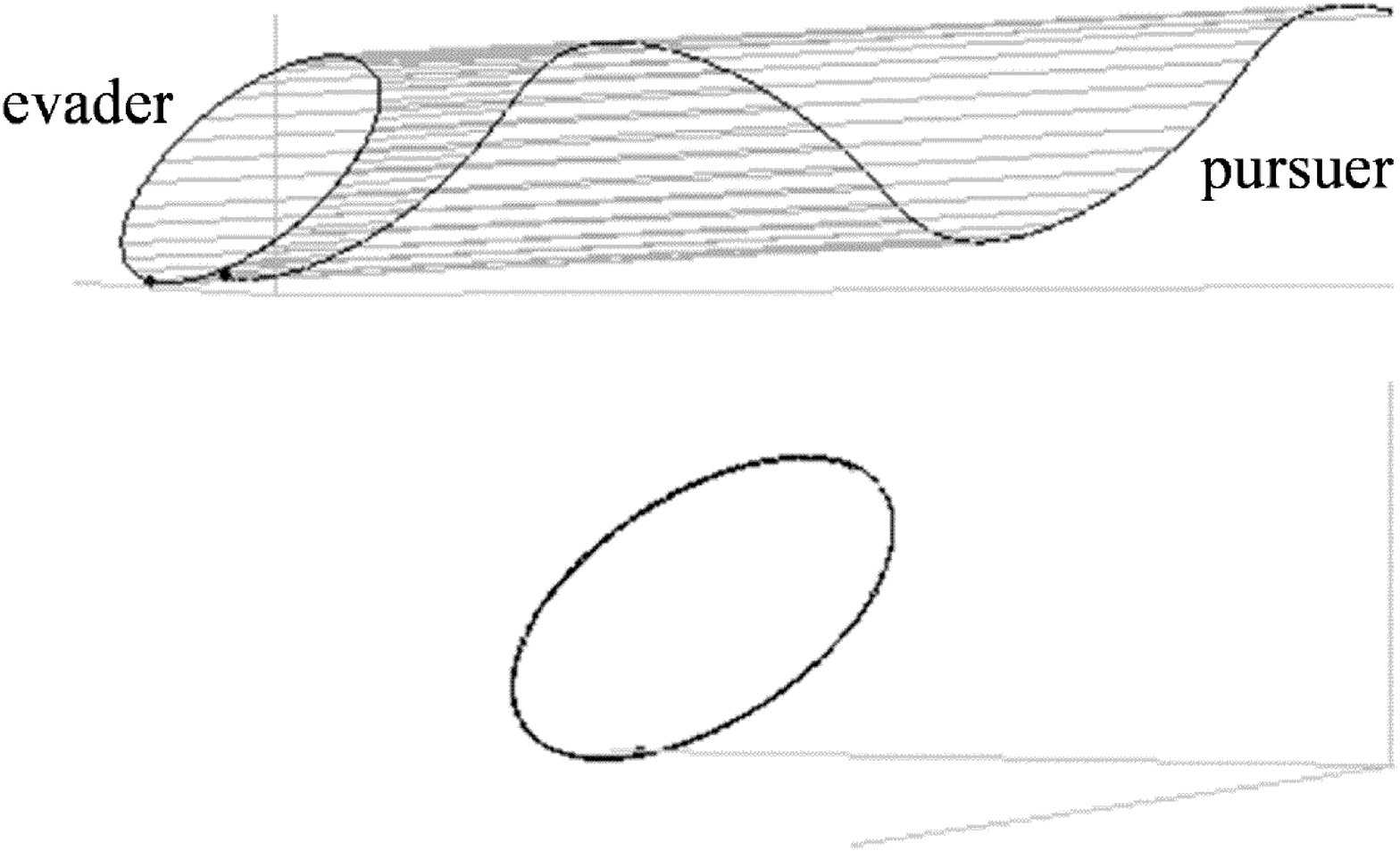}
\caption{\label{simfig4} Evader trajectory with constant curvature inputs
(i.e., a circular trajectory), and corresponding
pursuer trajectory.}
\end{figure}

Figures \ref{simfig1}-\ref{simfig4} illustrate the behavior of the
three-dimensional motion camouflage system (\ref{pursuer3d}), (\ref{evader3d})
under control law (\ref{adefn}) - (\ref{vpdefn})
for the pursuer, and various open-loop curvature controls for the
evader.  The speeds of the pursuer and evader are constant, and
the ratio of speeds is $ \nu_e/\nu_p = .9 $.  For each simulation,
two views of the resulting three-dimensional trajectories are shown:
one perpendicular to the ${\bf r}_\infty$-direction (upper
plot), and one along the $ {\bf r}_\infty$-direction (lower plot).  In 
figure \ref{simfig1}, the evader moves in a straight line (i.e., its
curvature controls are identically zero).  The corresponding motion
camouflage trajectory for the pursuer is then also a straight line.
The upper plot of figure \ref{simfig1} shows these straight-line
trajectories, along with the baselines at equally-spaced intervals of
time.  Recall that by definition, these baselines are parallel when
the system is in a state of motion camouflage.  In the lower plot of figure
\ref{simfig1}, the trajectories of the pursuer and evader overlap,
and the baselines are essentially normal to the page. 
 
In figure \ref{simfig2}, the curvature controls for the evader are
sinusoidal functions of time.  Whereas in figure \ref{simfig1}, the
motion is very nearly planar (with the plane determined by the initial
heading of the evader), in figure \ref{simfig2}, the motion is seen
to be truely three-dimensional.  Nevertheless, the baselines are
observed to be nearly parallel.
In figure \ref{simfig3}, the curvature controls for the evader are
randomly varying, and similarly to figure \ref{simfig2},
the trajectories are truly three-dimensional in character,
with the baselines nearly parallel.
In figure \ref{simfig4}, the curvature controls for the evader are
constant and nonzero, so that the trajectory of the evader is
circular.   

Although there is a brief transient period at the start of each
simulation during which $ \Gamma $ is driven close to $ -1 $ by the
control law, this transient period is such a small fraction of the
total simulation time that the transient behavior is not evident
in figures \ref{simfig1}-\ref{simfig4}.  The effect of the
gain $ \mu $ on both the duration of the transient and the
ultimate tolerance within which $ \Gamma $ remains near $ -1 $
is illustrated for the planar setting in \cite{mcarxiv05}.
Since the bounds and estimates for the three-dimensional problem
are analogous to the planar problem, similar behavior is expected.

\section{Connection to missile guidance}

For the planar setting, the connection between motion camouflage
and the {\it pure proportional navigation guidance} (PPNG) law
has been described in \cite{mcarxiv05}.  There is also a 
three-dimensional version of the PPNG law, which has been 
studied in \cite{songha} and \cite{oh_ha99}.
The PPNG law (by definition)
produces an acceleration which is perpendicular to the velocity of
the missile and proportional to the angular velocity
of the {\it line of sight} (LOS) vector.
If $ A_M $ denotes the lateral acceleration of the missile, 
$ V_M $ its velocity, and $ \Omega_L $ the angular velocity of 
the LOS vector, then the three-dimensional 
PPNG law is given by 
\begin{equation}
A_M^{\mathit PPNG} = N \left(\Omega_L \times V_M \right),
\end{equation}
where $ N > 0 $ is a dimensionless constant known as the navigation
constant \cite{songha}.

On the other hand, from equations
(\ref{omegaeqn}), (\ref{adefn}), (\ref{updefn}), and (\ref{vpdefn}),
we observe that for the motion camouflage law,
the lateral acceleration of the pursuer is
\begin{eqnarray}
A_M^{\mathit MCPG} \hspace{-.2cm} & = & \hspace{-.2cm}
  \mu \nu_p^2 {\bf a} =  \mu \nu_p^2 \left[ {\bf x}_p \times 
 \left(\dot{\bf r} \times \frac{\bf r}{|{\bf r}|}\right) \right] \nonumber \\
\hspace{-.2cm} & = & \hspace{-.2cm}
 - \mu \nu_p^2  |{\bf r}| 
 \left({\bf x}_p \times \mbox{\boldmath$\omega$\unboldmath} \right).
\end{eqnarray}
Identifying $ \Omega_L $ with $ \mbox{\boldmath$\omega$\unboldmath} $
and $ V_M $ with $ \nu_p {\bf x}_p $, we see that
\begin{equation}
A_M^{\mathit MCPG} = (\mu \nu_p |{\bf r}|) \left( \Omega_L \times V_M \right).
\end{equation}
To compare PPNG to MCPG, following the approach taken in the planar
setting \cite{mcarxiv05}, we take $ r_o $ to be a length scale for
the MCPG problem, and define the dimensionless gain
\begin{equation}
N^{\mathit MCPG} = \mu \nu_p r_o.
\end{equation}
Then
\begin{equation}
A_M^{\mathit MCPG} =
 \left(\frac{N^{\mathit MCPG} |{\bf r}|/{r_o}}{N} \right)A_M^{\mathit PPNG}.
\end{equation}
Thus, the MCPG law uses range information to provide high gain 
during the initial phase of the engagement, and ramps the gain
down to a lower value in the terminal phase ($|{\bf r}| \approx r_o $). 
This type of gain control is plausible for echolocating bats
(see \cite{ghose05}) which have remarkable ranging ability.

\vspace{-.05cm}

\section{Directions for futher work}

\vspace{-.05cm}

In the biological context, one direction being pursued is
the interpretation of three-dimensional trajectory data taken
from experiments in which a bat, {\it Eptesicus fuscus},
pursues a flying praying mantis (whose hearing organ
is disabled so that its trajectory is not influenced by the
presence of the bat).  The hypothesis that the bat uses an
MCPG strategy during the capture phase of its engagement with
the mantis is currently being tested using experimental
data collected in the Auditory Neuroethology Laboratory
at the University of Maryland (http://www.bsos.umd.edu/psyc/batlab).
This work represents part of a larger program to understand
sensory-motor processing and feedback in biological model systems. 

Another aspect of motion camouflage currently under study is
discovering feedback laws for motion camouflage with respect
to a finite point (as opposed to the point at infinity).
In finite-point motion camouflage, the pursuer uses a fixed
object as camouflage as it approaches the evader, and
this strategy also appears to be biologically revelant.
Various scenarios for motion camouflage involving teams of 
pursuers are also of interest, particularly in combination
with formation-control
laws based on gyroscopic interactions \cite{cdc05}. 
Some possible scenarios for team motion camouflage appear in
\cite{andersonteam}.


\vspace{-.15cm}

\begin{thebibliography}{99}

\vspace{-.1cm}

\bibitem{andersonneural} 
A.J. Anderson and P.W. McOwan, ``Model of a predatory
stealth behavior camouflaging motion,'' {\it Proc. Roy. Soc. Lond. B}
Vol. 270, No. 1514, pp. 489-495, 2003.

\bibitem{andersonteam} A.J. Anderson and P.W. McOwan, ``Motion camouflage
team tactics,''  Evolvability \& Interaction Symposium
(see http://www.dcs.qmul.ac.uk/\~{}aja/TEAM\_MC/team\_mot\_cam.html), 2003.


\bibitem{bishop} R.L. Bishop, ``There is more than one way to frame a curve,''
{\it The American Mathematical Monthly}, Vol. {82}, No. 3, pp. 246-251, 1975.

\bibitem{collett75} T.S. Collett and M.F. Land, ``Visual control of flight
behaviour in the hoverfly, {\it Syritta pipiens},'' {\it J. comp. Physiol.},
vol. 99, pp. 1-66, 1975.

\bibitem{ghose05} K. Ghose, T. Horiuchi, P.S. Krishnaprasad and C. Moss,
``Echolocating bats use a nearly time-optimal
strategy to intercept prey,'' {\it PLoS Biology}, to appear, 2006.

\bibitem{glendinning} P. Glendinning, ``The mathematics of motion
camouflage,''  {\it Proc. Roy. Soc. Lond. B}, Vol. 271, No. 1538, pp. 477-481,
2004.

\bibitem{cdc05} E.W. Justh and P.S. Krishnaprasad, ``Natural frames and
interacting particles in three dimensions,''   
{\it Proc. 44th IEEE Conf. Decision and Control}, 2841-2846, 2005
(see also arXiv:math.OC/0503390v1).

\bibitem{mcarxiv05} E.W. Justh and P.S. Krishnaprasad, ``Steering laws 
for motion camouflage,'' preprint, 2005 (arXiv:math.OC/0508023).

\bibitem{srini03} A.K. Mizutani, J.S. Chahl, and M.V. Srinivasan,
``Motion camouflage in dragonflies,'' {\it Nature}, Vol. 423, p. 604,
2003.

\bibitem{oh_ha99} J.H. Oh and I.J. Ha, ``Capturability of the 3-dimensional
pure PNG law,'' {\it IEEE Trans. Aerospace. Electr. Syst.}, vol. 35, No. 2,
pp. 491-503, 1999.

\bibitem{shneydor98} N.A. Shneydor, {\it Missile Guidance and Pursuit}, 
Horwood, Chichester, 1998.

\bibitem{songha} S.H. Song and I.J. Ha,``A Lyapunov-like approach to 
performance analysis of 3-dimensional pure PNG laws,'' {\it IEEE Trans.
Aerospace and Electronic Systems}, Vol. 30, pp. 349-358, 1994.

\bibitem{srini95} M.V. Srinivasan and M. Davey, ``Strategies for active
camouflage of motion,'' {\it Proc. Roy. Soc. Lond. B}, Vol. 259, No. 1354,
pp. 19-25, 1995.

\bibitem{srini04} M.V. Srinivasan and S. Zhang, ``Visual Motor
Computations in Insects,'' {\it Ann. Rev. Neurosci.}, Vol. 27, pp. 679-696,
2004.


\end{thebibliography}
\end{document}